\documentclass[a4paper,10pt,leqno,twoside]{article}

\usepackage[english]{babel}
\usepackage{german, umlaute, amsmath, amssymb , latexsym, 
  epic, epsfig, rotating,fancyheadings, amsthm}

\newcommand{\ld}{\ensuremath{,\ldots,}}
\newcommand{\ssq}{\ensuremath{\subseteq}}
\newcommand{\smin}{\ensuremath{\setminus}}

\newcommand{\ra}{\ensuremath{\rightarrow}}
\newcommand{\equi}{\ensuremath{\Leftrightarrow}}

\newcommand{\N}{\ensuremath{\mathbb{N}}} 
\newcommand{\R}{\ensuremath{\mathbb{R}}}
\newcommand{\Z}{\ensuremath{\mathbb{Z}}}
\newcommand{\Q}{\ensuremath{\mathbb{Q}}}

\newcommand{\kreis}{\ensuremath{\mathbb{T}^{1}}}
\newcommand{\ntorus}[1][2]{\ensuremath{\mathbb{T}^{#1}}}

\newcommand{\nLim}{\ensuremath{\lim_{n\rightarrow\infty}}}

\newcommand{\nfolge}[1]{\ensuremath{(#1)_{n\in\mathbb{N}}}}

\newcommand{\ind}{\ensuremath{\mathbf{1}}}

\newcommand{\inergsum}{\ensuremath{\sum_{i=0}^{n-1}}}

\newcommand{\ntel}{\ensuremath{\frac{1}{n}}}

\newcommand{\halb}{\ensuremath{\frac{1}{2}}}

\newcommand{\thx}{\ensuremath{(\theta,x)}}
\newcommand{\thom}{\ensuremath{\theta + \omega}}

\newtheorem*{thma}{Theorem A}
\newtheorem*{thmb}{Theorem B}
\newtheorem*{cora}{Corollary A}
\newtheorem{definition}{Definition}[section]
\newtheorem{thm}{Theorem} 
 
\newtheorem{lem}[definition]{Lemma}

\newtheorem{prop}[definition]{Proposition}

\newtheorem{bem}[definition]{Remark}

\numberwithin{equation}{section}

\title{\Large\textsc{Rotation numbers for quasiperiodically forced
    circle maps -- Mode-locking vs strict monotonicity}}
\author{Kristian Bjerkl\"ov\thanks{Department of Mathematics,
    University of Toronto,  Canada. Email:{\tt  bjerklov@math.utoronto.ca}} and
  Tobias J\"ager\thanks{Mathematisches Institut, Universit\"at
    Erlangen-N\"urnberg, Germany. Email: {\tt
      jaeger@mi.uni-erlangen.de}}} 

\newcommand{\Fhom}{\ensuremath{{\cal F}}}
\newcommand{\eps}{\ensuremath{\varepsilon}}

\begin{document}
\setlength{\topmargin}{-1cm}
\setlength{\oddsidemargin}{0.02\textwidth}  
\setlength{\evensidemargin}{0.02\textwidth}  
\setlength{\textheight}{0.9\textheight}

\maketitle \abstract{We describe the relation between the dynamical
  properties of a quasiperiodically forced orientation-preserving
  circle homeomorphism $f$ and the behavior of the fibered rotation
  number w.r.t.\ strictly monotone perturbations. Despite the fact
  that the dynamics in the forced case can be considerably more
  complicated, the result we obtain is in perfect analogy with the
  one-dimensional situation. In particular, the fibered rotation
  number behaves strictly monotonically whenever the rotation vector
  of $f$ is irrational.  This answers a question posed by Herman 
  \cite{herman:1983}.}

\section{Introduction}

For an orientation-preserving circle homeomorphism $g :\kreis \ra
\kreis$ it is well-known that the rotation number behaves strictly
monotonically whenever $g$ has no periodic points. On the other hand,
the only situation where mode-locking occurs (i.e.\ the rotation
number is stable w.r.t.\ perturbations) is when there exists a closed
interval which is mapped inside its own interior by some iterate of
$g$, as in the case of a stable periodic orbit. If the rotation number
only stays constant on one side this corresponds to the existence of
parabolic periodic points. In this paper we show that exactly
the same picture holds for quasiperiodically forced (qpf)
orientation-preserving circle homeomorphisms, where $p,q$-invariant
strips (as introduced in \cite{jaeger/keller:2006} and
\cite{jaeger/stark:2006}) serve as natural analogues of periodic
orbits. \bigskip

We consider continuous maps $f : \ntorus \ra \ntorus$ of 
the form
\begin{equation} \label{eq:qpfs}
  f\thx \ = \ (\thom,f_\theta(x)) \ ,
\end{equation}
where $\omega \in [0,1] \setminus \Q$. In addition, we require all
\textit{fiber maps} $f_\theta$ to be orientation-preserving circle
homeomorphisms. The class of all such maps will be denoted by $\Fhom$.
Further, we will use the notation $f^n_\theta(x) = \pi_2\circ
f^n\thx$.  Given any continuous lift $F :\kreis \times \R \ra \kreis
\times \R$ of $f \in \Fhom$, the limit
\begin{equation}
    \label{eq:rotnum}
        \rho(F) := \nLim (F_\theta^n(x) - x)/n
\end{equation}
exists and is independent of $\theta$ and $x$ \cite{herman:1983}.  The
\textit{fibered rotation number} of $f$ is defined as $\rho(f) :=
\rho(F) \bmod 1$. In order to state our results, we need the following
notions:
\begin{definition} \label{def:strictmonotonicity}
  Suppose $f$ is a qpf circle homeomorphism with lift $F$. We let 
  \[
      F_\eps(\theta,x) \ := \ (\thom,F_\theta(x)+\eps)\ 
  \]
  and say the rotation number is strictly monotone in $f$ if the map
  $\eps \mapsto \rho(F_\eps)$ is strictly monotone at $\eps = 0$.  If
  $\eps \mapsto \rho(F_\eps)$ is constant in a neighborhood of
  $\eps=0$ we say $f$ is mode-locked.
\end{definition}
\begin{thma} 
  Suppose $f\in {\cal F}$ has no $p,q$-invariant strip.  Then the
  rotation number is strictly monotone in $f$. In particular, this
  holds whenever $\omega$, $\rho(f)$ and 1 are rationally
  independent.%
\footnote{The presence of $p,q$-invariant graphs forces $\omega$,
$\rho(f)$ and 1 to be rationally dependent, see Theorem~\ref{thm:poincare}.
Hence the second statement is an immediate consequence of the first.}
\end{thma} \enlargethispage*{1000pt}
Note that the second statement in Theorem A is the answer to the
question asked by Herman in \cite[Section 5.18, p.\ 497]{herman:1983}.
\begin{thmb}
  Suppose $f \in {\cal F}$ is mode-locked. Then there exists a closed
  annulus, bounded by continuous curves, which is mapped into its own
  interior by some iterate of $f$.
\end{thmb} 
For a description of the critical cases with one-sided monotonicity, 
we refer to Section~\ref{ProofB}. \\

\pagebreak
\setlength{\topmargin}{+0.7cm}

There is one well-known situation where a similar result to that in
Theorem A holds, namely in the case of the (generalized) Harper map
\begin{equation}
  \label{eq:harper}
  s_E : \kreis \times \overline{\R} \ra \kreis \times \overline{\R}
  , \quad \thx \mapsto \left(\thom,V(\theta)-E-\frac{1}{x} \right) \ ,
\end{equation}
where $V:\kreis\to\mathbb{R}$ is a continuous function and
$E\in\mathbb{R}$. By letting $x_n=u_{n+1}/u_n$, there is a 1-1
correspondence between orbits of this map and formal eigenfunctions of
the discrete Schr\"odinger operator $(H_\theta
u)_n=-(u_{n+1}+u_{n-1})+V(\theta+n\omega)u_n$. Note that, by
identifying $\overline{\R}$ with \kreis, (\ref{eq:harper}) defines an
element of ${\cal F}$ and we can therefore speak of the fibered
rotation number of $s_E$. Due to a fruitful interplay between spectral
theory and dynamical systems methods it is well-known that the
function $E\mapsto \rho(S_E)$ is monotone and continuous, and maps
$\mathbb{R}$ onto $[0,1]$. In fact, it can be shown that
$\rho(s_E)=k(E)$, where $k(E)$ denotes the integrated density of
states for the corresponding Schr\"odinger operator.  Moreover, if
$\rho(s_E)$ is constant on an interval, then $\rho(s_E)=k\omega \bmod
1$ for some integer $k$ \cite{DS,johnson/moser:1982}.  In particular,
if $\rho(s_E)$ and $\omega$ are rationally independent, then
$\rho(s_E)$ is strictly monotone in $E$.  We would like to stress that
analogue results first where established for the continuous
Schr\"odinger equation with almost periodic potentials in the
fundamental paper \cite{johnson/moser:1982}.

In contrast to the situation for the Harper map and despite a large
number of numerical studies (\cite{prasad/negi/ramaswamy:2001} gives a
good overview and further reference), there are hardly any rigorous
results about other parameter families of qpf circle maps. The most
prominent example is probably the qpf Arnold circle map
\[
   f_{\alpha,\beta,\tau} : \ntorus \mapsto \ntorus \ , \quad \thx
   \mapsto \left(\thom,x+\tau+\frac{\alpha}{2\pi} \sin(2\pi x) +
   \beta g(\theta) \bmod 1\right)
\]
with parameters $\alpha,\tau \in [0,1],\ \beta \in \R$ and continuous
forcing function $g: \kreis \ra \R$. This map was first studied in
\cite{ding/grebogi/ott:1989}, where it was proposed as a simple model
of an oscillator forced at two incommensurate frequencies. It was
observed that the fibered rotation number seems to stay constant on
open regions in the $(\alpha,\tau)$-parameter space, the so-called
Arnold tongues (see also \cite{feudel/kurths/pikovsky:1995,
  glendinning/feudel/pikovsky/stark:2000,
  stark/feudel/glendinning/pikovsky:2002}).%
\footnote{For the unforced Arnold circle map this is certainly
  well-known.}
Although this was clearly expected and backed up by numerical
evidence, to the knowledge of the authors it was not known so far that
the rotation numbers have to be rationally related inside these
tongues. Now this is certainly a direct consequence of Theorem A,
which more generally implies the following
\begin{cora}
  Suppose $(f_\tau)_{\tau\in\R}$ is a parameter family in ${\cal F}$ with
  lifts $(F_\tau)_{\tau \in \R}$. Further assume $F_\tau$ depends
  continuously on $\tau$ and $F_{\tau,\theta}(x) <
  F_{\tau',\theta}(x) \ \forall \thx$ whenever $\tau < \tau'$. Then
  $\tau \mapsto \rho(f_\tau)$ can only stay constant in $\tau_0$ if
  $\omega$, $\rho(f_{\tau_0})$ and 1 are rationally dependent. 
\end{cora}
Finally, we note that Theorem A can also give interesting information
in situations with rationally related rotation numbers: The
simulations in \cite{feudel/kurths/pikovsky:1995,
  glendinning/feudel/pikovsky/stark:2000} suggest that for certain
parameters $\alpha,\beta$ the $\tau$-interval with fibered rotation
number zero is collapsed to a single point. This observation will be
confirmed in a forthcoming paper \cite{bjerkloev/jaeger:2006b}, by
showing that for suitable forcing function $g$ and parameters
$\alpha,\beta$ the map $f_{\alpha,\beta,0}$ has fibered rotation
number zero, but minimal dynamics and therefore no $p,q$-invariant
strips.  Consequently $\tau \mapsto \rho(f_{\alpha,\beta,\tau})$ is
strictly monotone in $\tau=0$ by Theorem A.

\textbf{Acknowledgments.} This cooperation started during a workshop
on `Dynamics of cocycles and one-dimensional spectral theory' at
Oberwohlfach 2005. We would like to thank the organizers David
Damanik, Russell Johnson and Daniel Lenz for making this possible.
K.B.\ was partially supported by SVeFUM.  T.J.\
was supported by grant Ke 514/6-1 of the German Science Foundation
(DFG), which also funded a visit of K.B. to Erlangen.


\section{Preliminaries} \label{Preliminaries}

For any set $A \ssq \ntorus$ and $\theta \in \kreis$ we let $A_\theta
:= \{ x \in \kreis \mid \thx \in A\}$. In the simplest case $p=q=1$,
the definition of an invariant strip is as follows:
\begin{definition}
  \label{def:invstrips} Let $f \in {\cal F}$. $A \ssq \ntorus$ is
  called a $1,1$-invariant strip if it is compact, $f$-invariant and
  for all $\theta \in \kreis$ the set $A_\theta$ consists of exactly
  one nontrivial interval (i.e.\ $\neq \emptyset$ or $\kreis$).
\end{definition}
Note that in particular this definition includes continuous invariant
curves. As the formulation for the general case $(p,q) \in \N^2$ is
slightly technical and we actually do not have to use it, we refrain
from giving the precise definition here and refer to
\cite{jaeger/keller:2006} or \cite{jaeger/stark:2006}. However, in
order to get a basic idea one should think of $p$ $p$-periodic closed
curves which are permuted by the action of $f$ and all wind around the
torus $q$ times in the $\theta$-direction. The union of these curves
intersects every fiber in exactly in $pq$ points, and roughly spoken
one allows each of these points to be replaced by a closed interval in
the definition of a $p,q$-invariant strip.

The important thing we have to know here is the fact that the
existence of a $p,q$-invariant strip $A$ forces the fibered rotation
number to be of the form
\begin{equation} \label{eq:rational}
  \rho(f) \ = \ \frac{k}{q} \omega + \frac{l}{pq} \ \bmod 1 \ ,
\end{equation}
where the integers $k,l,p$ and $q$ are determined by the topological
and dynamical structure of $A$ and vice versa (see Lemma~3.9 in 
\cite{jaeger/keller:2006}).  Further, by going over to a suitable
iterate, lifting the system to a $q$-fold cover $(\R / q\Z) \times
\kreis$ and rescaling, one can always transform a system with a
$p,q$-invariant graph into one with a $1,1$-invariant graph (compare 
Lemma~2.15 in \cite{jaeger/stark:2006}; this is the reason why we do 
not need the general definition). For the latter (\ref{eq:rational})
implies that the fibered rotation number will be of the form %
$\rho(f) = k\omega \bmod 1$, and by conjugating with %
$h : \thx \mapsto (\theta,x-k\theta)$ we can finally assume that
$\rho(f) = 0$. \bigskip

The proof of Theorem A is based on a classification result for qpf
circle homeomorphisms obtained in \cite{jaeger/stark:2006}. In order
to state it we have to introduce the \textit{`deviations from the
  constant rotation'}, which are defined as
\begin{equation} \label{eq:deviations}
  D_F(n,\theta,x) \ := \ F_\theta^n(x) - x - n\rho(F) \ .
\end{equation}
In contrast to the unforced situation, these deviations do not have to
be uniformly bounded in $(n,\theta,x)$ in the forced case. However, if
the deviations are unbounded, they have to be unbounded on every
single orbit (see Theorem~1.3 in
\cite{stark/feudel/glendinning/pikovsky:2002}). This motivates the
following definition:
\begin{definition}
  $f\in{\cal F}$ is called $\rho$-bounded if the deviations are
  uniformly bounded and $\rho$-unbounded otherwise.
\end{definition}
Note that the boundedness of the deviations does not depend on the
choice of the lift $F$. We also remark that, again by Theorem~1.3 in
\cite{stark/feudel/glendinning/pikovsky:2002}, in the $\rho$-unbounded
case there always exist orbits for which the deviations are unbounded
from above (and similarly from below). It turns out that in the
$\rho$-bounded case, the natural analogue of the Poincar\'e
Classification Theorem (e.g.\ \cite{katok/hasselblatt:1997}) holds:
\begin{thm}[Theorems 3.1 and 4.1 in \cite{jaeger/stark:2006}]
  \label{thm:poincare} 
  If $f\in {\cal F}$ is $\rho$-bounded, then either there exists a
  $p,q$-invariant strip and $\rho(f)$, $\omega$ and 1 are rationally
  dependent or $f$ is semi-conjugate to the irrational torus
  translation $\thx \mapsto (\thom,x+\rho(f))$ by a semi-conjugacy $h$
  which is fiber-respecting (i.e.\ $\pi_1\circ h = \pi_1$).
    If $f\in {\cal F}$ is $\rho$-unbounded, then neither of these
  alternatives can occur and the map is always topologically
  transitive.
 \end{thm}
 Examples of systems with $\rho$-unbounded behavior can be found e.g.\ 
 in \cite{herman:1983}, similar examples with minimal dynamics in
 \cite{bjerkloev:2005a}. \bigskip

For the proof of Theorem~B we have to introduce a number of notions
and facts concerning qpf monotone maps.  We call $F : \kreis \times \R
\ra \kreis \times \R$ a \textit{qpf monotone map} if it is continuous,
has skew product structure as in (\ref{eq:qpfs}) and all fiber maps
$F_\theta :\R \ra \R$ are monotonically increasing. In particular,
this is true if $F$ is the lift of some $f\in{\cal F}$.  Similar to
Definition~\ref{def:invstrips} we define an ($F$-)\textit{invariant
  strip} as a compact $F$-invariant set which consists of exactly one
non-empty interval on every fiber.  The \textit{upper} and \textit{lower bounding graphs}
of an invariant strip $A$ are defined as
\begin{equation}
  \label{eq:boundinggraph}
  \varphi^+_A(\theta) \ := \ \sup A_\theta \quad \textrm{and} \quad
  \varphi^-_A(\theta) \ := \ \inf A_\theta ,
\end{equation}
respectively. Due to compactness, $\varphi^+_A$ will be upper
semi-continuous (u.s.c.) and $\varphi^-_A$ lower semi-continuous
(l.s.c.). Note that more generally we can apply definition
(\ref{eq:boundinggraph}) to any bounded $F$-invariant set $A$ in order
to obtain invariant (but not necessarily semi-continuous) bounding
graphs.  If $A$ and $B$ are two invariant strips, we use the notation
\[
A \preccurlyeq B \quad :\equi \quad \varphi^-_A \leq \varphi^-_B \ \ 
\textrm{ and } \ \  \varphi^+_A \leq \varphi^+_B
\]
and
\[
    A \prec B \quad :\equi \quad \varphi^+_A  < \varphi^-_B \ .
\]
Further, for any two graphs $\varphi,\psi : \kreis \ra \R$ with
$\varphi \leq \psi$ we let
\[ 
[\varphi,\psi] \ :=  \{ \thx \mid x \in [\varphi(\theta),\psi(\theta)]\} \ 
,
\] 
similarly for open and half-open intervals. For any graph $\varphi :
\kreis \ra \R$ we denote its point set by the corresponding capital
letter, i.e.\ $\Phi := \{(\theta,\varphi(\theta)) \mid \theta \in
\kreis\}$, likewise for curves $\gamma : I \ra \R$ which are only
defined on a subinterval $I\ssq \kreis$. Further, we let
\begin{equation}
  \label{eq:phiminus}
  \varphi^+ \ := \ \varphi^+_{\overline{\Phi}} \quad \textrm{and}
  \quad \varphi^- \ := \ \varphi^-_{\overline{\Phi}} \ .
\end{equation}
For simplicity, we denote $\varphi^{+-} = (\varphi^+)^-$,
$\varphi^{-+} = (\varphi^-)^+$. 

We call a $F$-invariant strip $A$ \textit{minimal} whenever it does
not strictly contain any smaller $F$-invariant strip. Minimality of
$A$ is equivalent to the fact that $\varphi^{+-}_A = \varphi^-_A$ and
$\varphi^{-+}_A = \varphi^+_A$ (see~\cite{jaeger/stark:2006}). The
following lemma describes a simple procedure to obtain minimal strips
from semi-continuous invariant graphs (see Lemma~2.5 together with
Definition~2.7 and Remark~2.8(a) in \cite{jaeger/stark:2006}):
\begin{lem}
  \label{lem:reflexivegraphs}
  Suppose $F$ is a qpf monotone map and $\varphi$ is an u.s.c.\ 
  invariant graph. Then $[\varphi^-,\varphi^{-+}]$ is a minimal
  invariant strip. Similarly, if $\varphi$ is l.s.c.\ then
  $[\varphi^{+-},\varphi^+]$ is a minimal invariant strip.
\end{lem}
One of the most important properties of a minimal strip $A$ is the
fact that it is \emph{pinched}, meaning that there is at least one
fiber which intersects $A$ only in one
single point.%
\footnote{Obviously the invariance implies that the set ${\cal P}_A
  :=\{\theta \in \kreis \mid \textrm{card}(A_\theta) = 1 \}$ is dense
  in \kreis\ in this case, and from Baire's Theorem it follows that
  ${\cal P}_A$ is even residual \cite{stark:2003}.}
This follows from
\begin{lem}[Theorem~4.5 in \cite{stark:2003}]
  \label{lem:pinched} Suppose $\varphi$ is an upper
  semi-continuous invariant graph of a qpf monotone map. Then
  $[\varphi^-,\varphi]$ is pinched. Similarly, if $\varphi$ is lower
  semi-continuous then $[\varphi,\varphi^+]$ is pinched.
\end{lem}


\section{Strict monotonicity: Proof of Theorem~A} \label{ProofA}

From now on, we fix $f\in{\cal F}$ with lift $F_0$ and let $F_\eps$ be
as is Definition~\ref{def:strictmonotonicity}.
Theorem~\ref{thm:poincare} allows to divide the problem into two
cases, namely the $\rho$-unbounded one and the one with a
semi-conjugacy to an irrational translation. We treat them separately
in the following two lemmas, starting with the $\rho$-unbounded case:
\begin{lem} \label{lem:monotonicity}
  Suppose $f$ is $\rho$-unbounded. Then the rotation number
  is strictly monotone in $f$.
\end{lem}
\proof\ We only have to show that for any $\eps > 0$ there holds
$\rho(F_\eps) > \rho(F_0)$, as the case where $\eps < 0$ is symmetric.
Hence, fix any $\eps > 0$. Due to uniform continuity, there exists
$\delta > 0$ such that $|\theta-\theta'| < \delta$ implies
$F_{0,\theta}(x) \leq F_{\eps,\theta'}(x) \ \forall x \in \R$. By
induction and using the monotonicity of the fiber maps we obtain 
\begin{equation} \label{eq:bound}
  F_{0,\theta}^n(x) \ \leq \
F_{\eps,\theta'}^n(x) \quad \forall \theta,\theta': |\theta-\theta'| < 
\delta \ \forall x\in\R,n\in\N \ .
\end{equation}
Fix $N\in\N$ such that $\{n\omega\}_{n=0}^N$ is $\delta$-dense in the
circle. As we mentioned above, in the $\rho$-unbounded case there
always exists an orbit $(\theta_0,x_0)$ on which the deviations are
unbounded from above. Let $\theta_k := \theta_0 + k\omega$. As all
maps $F_\theta^n$ are lifts of circle homeomorphisms, obviously there
holds
\begin{equation}
|D(n,\theta,x)-D(n,\theta,x')| \ \leq \ 1 \quad \forall 
n\in\N,\theta\in \kreis ,x,x'\in\R \ .
\end{equation} 
Further, there exists a constant $C > 0$ which satisfies
\begin{equation}\label{eq:comparedeviations}
  |D(n,\theta,x)-D(n+m,\theta,x)| \ \leq \ C \quad \forall
  n\in\N,\theta\in\kreis,x\in\R \textrm{ and } m=1\ld N \ .
\end{equation}
As the deviations of the orbit of $(\theta_0,x_0)$ are unbounded from
above, we can choose $M \in \N$ with $D(M,\theta_0,x_0) \geq
C+2$. Together with (\ref{eq:comparedeviations}) this implies
\begin{equation} \label{eq:dev}
  D(M,\theta_k,x) \ \geq \ 1 \quad \forall k=0\ld N, x\in\R \ .
\end{equation}
Now we proceed inductively to show that for any $\theta\in\kreis,x \in
\R$ there holds
\begin{equation}
  F^{nM}_{\eps,\theta}(x) - x \ \geq \ n(M\rho(F_0)+1) \quad \forall n 
  \geq 0 \ ,
\end{equation}
which immediately implies $\rho(F_\eps) \geq \rho(F_0)+\frac{1}{M}$.
For $n=0$ there is nothing to prove, so assume the statement holds for
some $n\geq 0$. Let $x' := F_{\eps,\theta}^{nM}(x)$ and choose $k \in \{0\ld
N\}$ such that $|\theta+nM\omega-\theta_k| < \delta$ (recall the choice of
$N$). Then 
\begin{eqnarray*}
  F^{(n+1)M}_{\eps,\theta}(x) - x & = &  F^M_{\eps,\theta+nM\omega}(x') -
  x' + x' - x  \\ & \stackrel{(\ref{eq:bound})}{\geq} &
  F^M_{0,\theta_k}(x')-x' + n(M\rho(F_0)+1)\\ &
  \stackrel{\mbox{}}{=} & M\rho(F_0) + 
  D(M,\theta_k,x') + n(M\rho(F_0)+1) \\ &
  \stackrel{(\ref{eq:dev})}{\geq} &  (n+1)(M\rho(F_0)+1) \ .
\end{eqnarray*}

\qed

\ \\
This only leaves the case where $f$ is semi-conjugate to an irrational
torus translation. In order to treat this, it is convenient to look at
invariant measures. Let $h : \thx \mapsto (\theta,h_\theta(x))$ be the
fiber-respecting semi-conjugacy from Theorem~\ref{thm:poincare} and
define a measure $\mu$ on \ntorus\ by $\mu(A) = \lambda^2(h(A))\ 
\forall A \in {\cal B}(\ntorus)$, where $\lambda^n$ denotes the
Lebesgue measure on \ntorus[n]. Then it is easy to see that $\mu$ is
an ergodic $f$-invariant probability measure. (In fact, it can be
deduced from Theorem~4.1 in \cite{furstenberg:1961} that $\mu$ is the
only $f$-invariant probability measure in this situation.)  From the
definition of $\mu$ it follows that
\begin{equation} \label{eq:fibermeasures}
  \mu(A) \ = \ \int_{\kreis} \mu_\theta(A_\theta) \ d\theta \ \ \
  \forall A \in    {\cal B}(\ntorus) \  
\end{equation}
where $A_\theta := \{ x \in \kreis \mid \thx \in A\}$ and the
so-called fiber measures (or conditional measures) are defined as
$\mu_\theta(B) = \lambda^1(h_\theta(B)) \ \forall B \in {\cal
  B}(\kreis)$. Obviously, by this definition the measures $\mu_\theta$
are all continuous (in the sense that $\mu_\theta(\{x\}) = 0 \ \forall
x \in \kreis$) and $f_\theta$ maps $\mu_\theta$ to $\mu_{\thom}$,
i.e.\ $\mu_{\thom}=f^*_\theta\mu_\theta$. (A general discussion of
fiber measures can for example be found in \cite{arnold:1998}.)
Therefore the strict monotonicity of the rotation number in the
semi-conjugated case is a consequence of the following lemma.

\begin{lem}
  Suppose $f$ has an ergodic invariant probability measure
  $\mu$ with continuous fiber measures $\mu_\theta$. Then the rotation
  number is strictly monotone in $f$. In particular, this is true
  whenever $f$ is semi-conjugate to an irrational translation of the
  torus.
\end{lem}
\proof Fix $\eps > 0$. We identify $\mu$ and $\mu_\theta$ with their
natural lifts to $\kreis \times \R$ and $\R$, respectively. As the
fiberwise rotation number does not depend on \thx, it suffices to show
that there exist one $\thx \in \kreis \times \R$ such that
\begin{equation} \label{eq:rotnuminequality} 
  \nLim \ntel F^n_{\eps,\theta}(x) \ > \ \nLim \ntel F^n_{0,\theta}(x) \ .
\end{equation}
Let us first see that (\ref{eq:rotnuminequality}) is a consequence of the
following statement: 
\begin{equation} \label{eq:ptheta}
  p(\theta) := \min\{ p \in\N \mid F^p_{\eps,\theta}(x)  \geq  F_{0,\theta}^p(x)
  + 1 \ \forall x \in \R \} \ < \ \infty \ \  \textit{for }
  \lambda^1\textit{-a.e. } \theta \in \kreis \ .
\end{equation}
Indeed, this implies that for some $q \in \N$ the set $A_q := \{
\theta \in \kreis \mid p(\theta)=q \}$ has positive measure, as
$\lambda^1\left(\bigcup_{q\in\N} A_q\right) = 1$. W.l.o.g.\ we can
assume that $q=1$, otherwise we replace $f$ by its $q$th iterate. Due
to the monotonicity and periodicity of the fiber maps we obtain
\[
F_{\eps,\theta}^{n}(x) \ \geq \ F_{0,\theta}^{n}(x) + \inergsum
\ind_{A_1}(\theta+i\omega) \ .
\]
As $\theta \mapsto \thom$ is ergodic, this implies
\[ 
 \nLim \ntel F_{\eps,\theta}(x) \ \geq \ \nLim \ntel F^n_{0,\theta}(x) +
 \lambda^1(A_1) \ > \ \nLim \ntel F^n_{0,\theta}(x)
\]
for $\lambda^1$-a.e.\ $\theta$, thus proving (\ref{eq:rotnuminequality}).

\ \\ It remains to show that the function $p$ is $\lambda^1$-a.s.\ finite. To
that end, note that for any $\eps>0$ the function
\[
    g\thx \ := \ \mu_\theta([x,x+\eps))
\]
is $\mu$-a.s.\ strictly positive.  Therefore
\[
\delta \ := \ \halb \int_{\ntorus} g \ d\mu \ > \ 0 \ .
\]
By Birkhoff's Ergodic Theorem, for $\mu$-a.e.\ $\thx$ we have
\begin{equation} \label{eq:Fergsum}
   \nLim \ntel \inergsum g\circ f^i\thx \ = \ 2\delta 
    \ ,
\end{equation}
and from (\ref{eq:fibermeasures}) it follows that (\ref{eq:Fergsum})
holds $\mu_\theta$-a.s.\ for $\lambda^1$-a.e.\ fixed $\theta$.
Let $S \ssq \kreis$ be a set of measure $\lambda^1(S)=1$ which is invariant
under rotation by $\omega$ and has the property that
(\ref{eq:Fergsum}) holds for $\mu_\theta$-a.e.\ $x\in \kreis$ whenever
$\theta \in S$.  Then for every $\theta \in S$ and $\mu_\theta$-a.e.
$x$ there holds
\[
\kappa\thx \ := \ \min\{ k\in\N \mid g\circ f^k\thx > \delta\} 
 \ < \ \infty \ .
\]
We will show that $p$ only takes finite values on $S$. In order to do so, we
fix $\theta_0 \in S$ and $x_0 \in \R$ such that $\kappa(\theta_0,x_0)
< \infty$.%
\footnote{Note that we apply $g$ and $\kappa$ also to points $\thx \in
  \kreis \times \R$ by identifying $x \in \R$ with its
  projection to \kreis.}
We will construct sequences $k_i$ and $y_i$ ($i\in\N_0$) with the
following properties: If $K_i := \sum_{j=0}^i k_j$, $\theta_i :=
\theta_0 + K_i\omega,\ x_i := F_{\theta_0}^{K_i}(x_0)$ and $z_i :=
F_{\eps,\theta_0}^{K_i}(x_0)$ then
\begin{eqnarray}
  x_i \ \leq \ y_i & \leq & z_i \\
  \mu_{\theta_i}([x_i,y_i)) & \geq & i \cdot \delta
  \label{eq:measuregrowth} \\
   \kappa(\theta_i,y_i) & < & \infty \ . \label{eq:kappafinite}
\end{eqnarray}
If $i > 2/\delta$ then (\ref{eq:measuregrowth}) gives
\[
  \mu_{\theta_i}([x_i,z_i)) \ \geq \ \mu_{\theta_i}([x_i,y_i)) \
  > \ 2
\] and consequently, as $\mu_\theta([x,x+1)) = 1 \ \forall \thx$
(recall that $\mu$ is a probability measure)
\[
  F_{\eps,\theta_0}^{K_i}(x_0) \ > \ F_{\theta_0}^{K_i}(x_0) + 2 \ ,
\]
which in turn implies $p(\theta_0) \leq K_i$. As $\theta_0 \in S$ was
arbitrary and $\lambda^1(S)=1$, this shows (\ref{eq:ptheta}).

The sequences $k_i, y_i$ are constructed as follows: If we let
$y_0=x_0$ and $k_0=0$, then for $i=0$, there is nothing to show.
Suppose $k_0 \ld k_i$ and $y_0 \ld y_i$ with the required properties
have been chosen. Let $k_{i+1} := \kappa(\theta_i,y_i)$. Then
\begin{eqnarray*}
    \lefteqn{\mu_{\theta_{i+1}}([x_{i+1},z_{i+1})) \ \geq \
    \mu_{\theta_{i+1}}([x_{i+1},F_{\eps,\theta_i}^{k_{i+1}}(y_i)))} \\
    & = & \mu_{\theta_i}([x_i,y_i)) +
    \mu_{\theta_{i+1}}([F_{\theta_i}^{k_{i+1}}(y_i),
    F_{\eps,\theta_i}^{k_{i+1}}(y_i) )) \ > \ i\delta + \delta \ .
\end{eqnarray*}
Thus we can choose some $y_{i+1} \in (x_{i+1},z_{i+1})$ which
satisfies (\ref{eq:measuregrowth}) and (\ref{eq:kappafinite}).

\qed


\section{Mode-locking: Proof of Theorem B} \label{ProofB}

For $f\in{\cal F}$ with lift $F_0$ let $F_\eps$ be as is
Definition~\ref{def:strictmonotonicity}. Due to Theorem A, we only
have to consider the case where $f$ has a $p,q$-invariant strip.
Further, as mentioned in Section~\ref{Preliminaries} we can always
assume that $p=q=1$ and $\rho(f)=0$. In this situation, the statement
of Theorem B is a consequence of the following:
\begin{prop}\label{prop:onesided}
  Suppose $\rho(f) = \rho(F_0) = 0$ and for some $\eps > 0$ there holds
  $\rho(F_\eps) = 0$. Then there exists a continuous curve $\gamma^+ :
  \kreis \ra \R$ which is mapped strictly below itself by some iterate
  of $F_0$. Similarly, if $\rho(F_{-\eps})=0$ then there exists a curve
  $\gamma^-$ which is mapped strictly above itself by some iterate of
  $F_0$.
\end{prop}
For if $f$ is mode-locked, then we obtain two curves $\gamma^+$ and
$\gamma^-$ which are mapped strictly above, respectively below
themselves by some iterate of $F_0$. It is easy to see that $\gamma^+$
and $\gamma^-$ must be disjoint and that they are lifts of the
continuous boundaries of a closed annulus ${\cal A} \ssq \ntorus$
which is mapped inside its own interior by some iterate of $F_0$.

\begin{bem} As we reduced the problem to the situation with fibered
  rotation number 0, the closed annulus from
  Proposition~\ref{prop:onesided} will be of `homotopy type' $(1,0)$,
  meaning that it only winds once around the torus, in the
  $\theta$-direction. In the general case, by redoing the
  transformations described in Section~\ref{Preliminaries}, we obtain
  a closed annulus of homotopy type $(q,k)$ with $q$ and $k$ as in
  (\ref{eq:rational}).
  
  Proposition~\ref{prop:onesided} also contains the description of the
  one-sided cases: Suppose the rotation number only remains constant
  on one side, for example $\rho(F_\eps) = 0$ for some $\eps >0$ but
  $\rho(F_{\eps'}) < 0 \ \forall \eps'<0$. Then the proposition still
  yields the existence of a closed curve $\gamma^+$ which is mapped
  below itself by some iterate of $F_0$.  The iterates of $\gamma^+$
  form a monotonically decreasing sequence of continuous curves. In
  the limit, they have to converge to an upper semi-continuous graph
  $\varphi^+$ which is the upper boundary of an invariant strip.
  Projecting the situation to the torus, we obtain the existence of a
  closed curve which is mapped towards some $1,1$-invariant strip in
  the clockwise direction.  (Note that `being mapped above or below
  itself' does not make any sense on the torus, unless there is some
  kind of reference object.) As there is no curve which is mapped
  towards this strip in the counterclockwise direction, the situation
  resembles the one with a parabolic periodic orbit in the
  one-dimensional case. Again, this picture remains valid in the case
  of a general $p,q$-invariant strip, only the homotopy type of all
  objects involved changes.
\end{bem}
Before we turn to the proof of Proposition~\ref{prop:onesided}, we
need two auxiliary lemmas. The first one concerns the possible
dynamics between two neighboring minimal invariant strips. As the
proof is basically the same as in \cite[Part 2 of the proof of Thm.
4.4]{jaeger/keller:2006}, we keep it rather brief and refer to
\cite{jaeger/keller:2006} for more details.
\begin{lem} \label{lem:B1}
  Suppose $A \prec B$ are two minimal $F_0$-invariant strips and that
  there is no invariant strip strictly in between. Let $C :=
  (\varphi^+_A,\varphi^-_B)$. Then either $F_{0|C}$ is topologically
  transitive or there exists a wandering open set $W \ssq C$.
\end{lem}
\proof Suppose neither of these two alternatives holds. Then, as $F_{0|C}$
is not transitive, there exist open balls $U,V$ such that the
invariant sets $\tilde{U}:= \bigcup_{n\in\Z} F_0^n(U)$ and $\tilde{V} :=
\bigcup_{n\in\Z} F_0^n(V)$ are disjoint. As $U$ is non-wandering we can assume,
by going over to a suitable iterate, that $F_0(U) \cap U \neq
\emptyset$. Thus, the successive iterates of $U$ are attached to each
other `like beads on a string' and form an open and path-connected
`tube' which winds around the torus. As any open subset of $U$ is also
non-wandering, eventually this tube has to close, meaning that some
higher iterate $F_0^n(U)$ intersects $U$ again after the tube has already
performed a whole loop.  Hence, the set $\tilde{U}$ has the property
that for each of its points it contains a closed curve going through
this point and winding at least once around the torus.

By going over to an even higher iterate if necessary, we can repeat
this argument and obtain the same property for $\tilde{V}$. But as
both $\tilde{U}$ and $\tilde{V}$ are connected, this implies that one
of the sets must lie strictly above the other, w.l.o.g.\ 
$\varphi^+_{\tilde{U}} \leq \varphi^-_{\tilde{V}}$.  Hence, as both
sets are invariant the set
$[\varphi^+_{\tilde{U}},\varphi^-_{\tilde{V}}]$ is an invariant strip,
contradicting the assumption that there is no such strip between $A$
and $B$.

\qed

The following lemma allows to draw further conclusions from the
existence of a wandering set:
\begin{lem} \label{lem:B2}
  Let $A$, $B$ and $C$ be as in Lemma~\ref{lem:B1} and suppose there
  is a wandering open set $W \ssq C$. Then there exists a curve $\gamma :
  \kreis \ra \R$ with $\Gamma := \{(\theta,\gamma(\theta)) \mid \theta
  \in \kreis\} \ssq C$ which is mapped either strictly above or
  strictly below itself by some iterate of $F_0$.
\end{lem}
\proof We construct the point set $\Gamma$ of the curve $\gamma$, as
this will greatly simplify the argument.  Note that we may allow
$\Gamma$ to contain vertical segments: The property we are interested
in is open w.r.t.\ Hausdorff distance. Therefore any vertical parts
can be slightly tilted in order to obtain a curve that can be
represented as a graph over $\kreis$.

By going over to a suitable iterate, we can assume w.l.o.g.\ that $W$
contains a ball of diameter larger than $\omega$. Let $\Lambda_0 \ssq
W$ be a straight horizontal line segment of length $\omega$.
Denote the endpoints of $\Lambda_0$ by $a_0=(\theta_0,x_0)$ and %
$b_0 = (\theta_0+\omega,x_0)$. Further, let %
$\Lambda_k := F_0^k(\Lambda_0)$, $a_k = (\theta_k,x_k) := F_0^k(a_0)$,
$b_k = (\theta_{k+1},y_k) := F_0^k(b_0)$ and denote the closed
vertical line segment between $b_k$ and $a_{k+1}$ by $[b_k,a_{k+1}]$,
similarly for open and half-open segments. Finally, we define %
$\Gamma_k$%
$:= \bigcup_{i=0}^k \Lambda_i \cup \bigcup_{i=0}^{k-1} [b_i,a_{i+1}]$.
Note that thus $\Gamma_k$ is a curve that joins $a_0$ and $b_k$ and
contains $k$ vertical segments.  Further, by construction
\begin{equation} \label{eq:gammaconstr1}
 F_0(\Gamma_k) \smin \Gamma_k \ =  \ \Lambda_{k+1} \cup (b_ka_{k+1}] \ .
\end{equation} 
Now let $n := \max\{k\in\N \mid k\omega < 1\}$ and define $\Gamma$ as
the union of $\Gamma_n \smin \pi_1^{-1}([0,\omega))$, $\Lambda_0$ and
the vertical line segment between $a_0$ and $\Lambda_n$. Note that
this defines a closed curve that winds exactly once around the torus
and contains $n+1$ vertical segments. Further, if we assume w.l.o.g.\ 
$n\geq 2$, then $\Gamma_1 \ssq \Gamma$. From now on we assume that
$a_1$ lies above $b_0$, such that $a_{k+1}$ lies above $b_k$ for all
$k\in\N_0$ by monotonicity. (If $a_1$ lies below $b_0$, this can be
treated similarly.) We have to distinguish three different cases:

First, assume there exists $m\in\N$ such that $b_m$ lies above
$\Lambda_0$.%
\footnote{Suppose $\Lambda$ is the point set of a curve defined on a
  subinterval of \kreis\ and $b\in\ntorus$. Then by saying
  $b=(\theta,x)$ is above $\Lambda$ we mean $x> \sup \Lambda_\theta$,
  implicitly assuming $\pi_1(b) \in \pi_1(\Lambda)$. In the following
  we will use similarly obvious terminology without further
  explanation.}
W.l.o.g.\ we can assume $m=n$, otherwise we lift the system to the
$j$-fold cover $(\R/j\Z) \times \kreis$ of $\ntorus$, where $j$ is the
integer part of $m\omega$, and repeat the construction of $\Gamma$ as
above. (It is obvious that if there exists a curve with the required
property on this $j$-fold cover, then the same is true for the
original system.)  If $b_n$ lies above $\Lambda_0$, by the assumption
made above the same is true for $a_{n+1}$.  Further, by the
monotonicity of the fiber maps $b_{n+1}$ lies above $\Lambda_1$. As
$\Lambda_{n+1}$ joins $a_{n+1}$ and $b_{n+1}$ and cannot intersect
$\Lambda_0 \cup \Lambda_1$ (recall that $\Lambda_0$ is contained in
the wandering set $W$), we obtain that $(b_n,a_{n+1}] \cup
\Lambda_{n+1}$ lies strictly above the curve $\Gamma_1$.  Together
with (\ref{eq:gammaconstr1}) this implies $F_0(\Gamma) \succcurlyeq
\Gamma$ and $F_0^N(\Gamma) \succ \Gamma$ for suitably large $N$.

Secondly, assume there exists $m\in\N$ such that $b_m$ lies below
$\Lambda_0$ and $\Lambda_{m+1}$ does not intersect $[b_0,a_1]$. Again,
we can assume $m=n$. As $\Lambda_{n+1}$ cannot intersect $\Lambda_0
\cup \Lambda_1$, it is disjoint from $\Gamma_1$ in this case. But as
$b_{n+1}$ lies below $\Lambda_1 \ssq \Gamma_1$ by monotonicity, this
implies that the whole curve $\Lambda_{n+1} \cup (b_n,a_{n+1}]$ is
below $\Gamma_1$.  Similar to above we obtain $F_0(\Gamma)
\preccurlyeq \Gamma$ and $F_0^N(\Gamma) \prec \Gamma$ for suitably
large $N$.

Finally, we have to address the case where $b_m$ lies below
$\Lambda_0$ and $\Lambda_{m+1}$ intersects $[b_0,a_1]$ whenever
$\pi_1(b_m) \in \pi_1(\Lambda_0)$. We show that this implies the
existence of an invariant strip strictly between $A$ and $B$,
contradicting the assumption that there is no such strip.

Let $\Omega := \bigcap_{k\geq0} \overline{\bigcup_{j=k}^\infty
  \Gamma_j}$.  Clearly $\Omega$ is a non-empty, compact,
$F_0$-invariant set, and consequently $\varphi^+_\Omega$ is an upper
semi-continuous invariant graph.  We claim that
\begin{equation}
  \label{eq:phiomega}
  x_0 + \delta \ \leq \ \varphi^+_\Omega(\theta) \ \leq \ x_1 + \delta 
  \ \ \ \forall \theta \in I := (\theta_1-\delta,\theta_1) \ 
\end{equation}
for a suitably small $\delta > 0$. This implies that the same
inequalities hold for $\varphi^{+-}_\Omega$, such that
$[\varphi^+_\Omega,\varphi^{+-}_\Omega]$ defines an invariant strip
which lies strictly between $A$ and $B$. (Note that whenever a strip
lies strictly between $A$ and $B$ on an open interval, then this is
true on all of \kreis.) 

As $\Gamma_0$ is contained in the wandering set $W$, there exist small
boxes $W_0 := B_\delta(\theta_1) \times B_\delta(x_0)$ around $b_0$
and $W_1 := B_\delta(\theta_1) \times B_\delta(x_1) $ around $a_1$
which no curve $\Lambda_j$ with $j \geq 2$ can intersect. We fix
$\delta \in (0,\omega)$ with this property. Now, whenever $\pi_1(b_m)
\in \pi_1(\Lambda_0)$ the curve $\Lambda_m \cup [b_m,a_{m+1}] \cup
\Lambda_{m+1} = F_0^m(\Gamma_1)$ has to pass through below the set
$\widehat{W}_1 := I \times B_\delta(x_1)$: This holds for
$\Lambda_{m+1}$ as this curve must intersect $[b_0,a_1]$ and cannot
intersect $\widehat{W}_1$, and $\Lambda_m$ lies below $\Lambda_0$
anyway as this is true for its right endpoint $b_m$. (Recall here that
we assumed that $b_k$ is always below $a_{k+1}$, such that $a_1$ is
above $\Lambda_0$.) Consequently none of the sets $\bigcup_{j\geq k}
\Gamma_j$ intersects the region $\{\thx \mid \theta \in I,\ x >
x_1-\delta\}$, and from this the upper bound in (\ref{eq:phiomega})
follows easily.

For the lower bound, note that there are infinitely many $m\in\N$ such
that $\pi_1(W_1) \ssq \pi_1(\Lambda_{m+1})$. For these,
$\Lambda_{m+1}$ has to pass through between the boxes $W_0$ and $W_1$
on their whole width.  Therefore the upper bounding graphs of the sets
$\bigcup_{j\geq k} \Gamma_j$ are always above $W_0$, and consequently
the same is true for their pointwise limit $\varphi^+_\Omega$.

\qed

\proof[Proof of Proposition~\ref{prop:onesided}] Fix $f\in{\cal F}$
with lift $F_0$ and suppose $\rho(F_0) = 0$ and $f$ has an invariant
strip. Further, assume there exists no continuous curve $\gamma^+ :
\kreis \ra \R$ which is mapped strictly below itself by some iterate
of $F_0$. We have to show that in this case there holds $\rho(F_\eps)
> 0$ for all $\eps > 0$. First of all, note that it is sufficient to
prove that
\begin{equation}
  \label{eq:prop1}
   \sup_{\theta,x,n} (F^n_{\eps,\theta}(x) - x) \ = \ \infty \quad \forall
   \eps > 0 \ ,
\end{equation}
where the supremum is taken over all $\thx \in \kreis \times \R$ and
$n\in\N$. Indeed, if $\eps'>0$ is fixed we can apply (\ref{eq:prop1})
to $\eps'/2$ and obtain two possibilities: First, we could have
$\rho(F_{\eps'/2}) > 0$, but in this case we are finished as
$\rho(F_{\eps'}) \geq \rho(F_{\eps'/2}) > 0$. The only other
alternative is to have $\rho(F_{\eps'/2}) = 0$ but unbounded
deviations for $F_{\eps'/2}$. However, in this case $\eps \mapsto
\rho(F_\eps)$ is strictly monotone in $\eps = \eps'/2$ due to
Lemma~\ref{lem:monotonicity}, such that again $\rho(F_{\eps'}) > 0$.

Fix $\eps>0$. In order to prove (\ref{eq:prop1}), we will show that an
orbit of $F_\eps$ moves upwards along a suitable sequence of minimal
$F_0$-invariant strips $A_n$. As we want to proceed by induction and
there might be uncountably many invariant strips, we have to make a
certain choice in the construction of this sequence. To that end,
given any minimal invariant strip $A$ we denote by ${\cal M}(A)$ the
set of all minimal invariant strips $B \succ A$. (Note that ${\cal
  M}(A)$ is always non-empty, as all integer translates of $A$ are
minimal invariant strips as well.) Further, if $B \in {\cal M}(A)$,
let
\[
    D(A,B) \ := \ \inf_{\theta \in \kreis} \left(\varphi^-_B(\theta) -
    \varphi^+_A(\theta) \right)
\]
and 
\[
{\cal M}_\eps(A) \ := \ \left\{ B \in {\cal M}(A) \mid D(A,B) \leq
\frac{\eps}{2} \right\} \ .
\]
We start the construction of the sequence \nfolge{A_n} with any
minimal invariant strip $A_0$. (Note that $1,1$-invariant strips of
$f$ lift to $F_0$-invariant strips, so by assumption such an $A_0$
always exists.) Then we proceed by induction in the following way:
Suppose $A_0 \prec A_1 \prec \ldots \prec A_{n-1}$ have been chosen.
We distinguish two cases:

\ \\
$({\cal A}1)$: If ${\cal M}_\eps(A_{n-1})$ is empty, we choose $A_n$
as the (unique) minimal invariant strip above $A_{n-1}$ with the
property that there is no other invariant strip strictly between
$A_{n-1}$ and $A_n$.%
\footnote{$A_n$ can be constructed as follows: Let
  $\varphi(\theta) := \inf_{B\in{\cal M}(A_{n-1})} \varphi^+_B(\theta)$.
  Then $\varphi$ is an upper semi-continuous invariant graph, and
  $A_n:=[\varphi^-,\varphi^{-+}]$ is a minimal invariant strip
  (see Lemma~\ref{lem:reflexivegraphs}) which has the
  required property.}

\ \\
$({\cal A}2)$: If ${\cal M}_\eps(A_{n-1})$ is non-empty we let
\begin{equation} \label{eq:phin+1}
\psi(\theta) \ := \ \sup_{B \in {\cal M}_\eps(A_{n-1})} \varphi^-_B(\theta) 
\end{equation}
and 
\begin{equation}  \label{eq:An+1}
   A_n \ := \ [\psi^{+-},\psi^+] \ .
\end{equation}
For sufficiently large $c\in\R$, $\psi$ is the lower bounding graph of 
the compact set 
\[
{\cal C} \ := \ \bigcap_{B \in {\cal M}_\eps(A_{n-1})} [\varphi^-_B,c] 
\ .
\]
Thus $\psi$ is l.s.c. and therefore $A_n$ is a minimal invariant strip
by Lemma~\ref{lem:reflexivegraphs}. Further, it follows from the
finite intersection property that ${\cal C}$ intersects the set
$[\varphi^-_{A_{n-1}},\varphi^+_{A_{n-1}}+\eps/2]$. Hence, there
exists at least one $\theta'\in \kreis$ with
\begin{equation} \label{eq:theta'}
\psi(\theta') - \varphi^+_{A_{n-1}}(\theta') \ \leq \ \eps/2 \ .
\end{equation}

\ \\ 
Note that constructing the sequence $A_n$ in this way we obtain 
$D(A_{n-1},A_{n+1}) \geq \eps/2$ for all $n \geq 1$. Thus,
\[
\lim_{n \ra \infty} \inf_{\theta\in\kreis} \varphi^+_{A_n}(\theta) \ =
\ \infty \ .
\] 
Therefore (\ref{eq:prop1}) is an immediate consequence of the
following claim, which concludes the proof of the proposition:

\ \\
\textbf{Claim:} \ For any $n \geq 1$, $\theta_0 \in \kreis$ and $x_0 >
\varphi^+_{A_{n-1}}(\theta_0)$ there exists $k \in \N$ which satisfies
\[
F^k_{\eps,\theta_0}(x_0) \  > \ \varphi^+_{A_n}(\theta_0+k\omega) \ .
\]
\proof[Proof of the Claim] First of all, we fix $\delta_1 > 0$ such
that
\begin{equation} \label{eq:bound2}
  F_{0,\theta}^n(x) \ + \ \eps/2 \leq \
F_{\eps,\theta'}^n(x) \quad \forall \theta,\theta': |\theta-\theta'| < 
\delta_1 \ \forall x\in\R,\ n \geq 1 \ .
\end{equation}
(Compare (\ref{eq:bound}).) Further, we have to distinguish three
cases. (See Lemmas~\ref{lem:B1} and \ref{lem:B2} for the division of
$({\cal A}1)$ into the following two sub-cases.)

\ \\
\underline{Case $({\cal A}1)(a)$:} \ \emph{There is no invariant strip
  strictly between $A_{n-1}$ and $A_n$ and $F_{0}$ restricted to
  $(\varphi_{A_{n-1}}^+,\varphi_{A_n}^-)$ is topologically
  transitive.}  \bigskip

Choose $\theta^*$ such that $A_n$ is pinched at $\theta^*$, i.e.\ 
$\varphi^-_{A_n}(\theta^*) = \varphi^+_{A_n}(\theta^*)$ (recall that
$A_n$ is pinched by Lemma~\ref{lem:pinched}).  As $A_n$ is compact we
can further choose $\delta_2 \in (0,\delta_1)$ such that
\begin{equation}
\label{eq:pinchedball}
  [\varphi^-_{A_n}(\theta),\varphi^+_{A_n}(\theta)] \ \ssq \
  B_{\eps/8}(\varphi^+_{A_n}(\theta^*)) \quad \forall \theta
  \in B_{\delta_2}(\theta^*) \ .
\end{equation}
Let $(\theta',x')$ be a point with dense orbit in
$(\varphi_{A_{n-1}}^+,\varphi_{A_n}^-)$, $\theta' \in
B_{\delta_2/2}(\theta_0)$ and $x' \leq x_0$. Choose $k \in \N$ such
that $\theta'+k\omega \in B_{\delta_2/2}(\theta^*)$ and
$F^k_{0,\theta'}(x') > \varphi^+_{A_n}(\theta^*) - \eps/4$. Then
  \begin{eqnarray*}
    F_{\eps,\theta_0}^k(x_0) & \geq & F_{\eps,\theta_0}^k(x') \
    \stackrel{(\ref{eq:bound2})}{\geq} F^k_{0,\theta'}(x') + \frac{\eps}{2} \
    > \ \varphi^+_{A_n}(\theta^*) + \frac{\eps}{4} \
    \stackrel{(\ref{eq:pinchedball})}{>}  \
    \varphi^+_{A_n}(\theta_0+k\omega) \ .
  \end{eqnarray*}

\ \\
\underline{Case $({\cal A}1)(b)$:} \ \emph{There is no invariant strip
strictly between $A_{n-1}$ and $A_n$, but there is a curve $\gamma :
\kreis \ra \R$ which is mapped strictly above itself by some
iterate of $F_0$.} \bigskip

Define $\gamma_n(\theta) :=
F^n_{\theta-n\omega}(\gamma(\theta-n\omega))$. As this sequence of
continuous curves is monotonically increasing and bounded, it must
converge pointwise to some l.s.c.\ invariant graph $\gamma_\infty :=
\lim_{n\ra\infty} \gamma_n$. Note that $\gamma_\infty$ does not necessarily
have to coincide with $\varphi^-_{A_n}$. However, we must have
$\gamma_\infty^+ = \varphi^+_{A_n}$, otherwise
$[\gamma_\infty,\gamma_\infty^+]$ would be an invariant strip strictly 
between $A_{n-1}$ and $A_n$. Therefore, by
Lemma~\ref{lem:pinched} the set %
$\tilde{A} := [\gamma_\infty,\varphi^+_{A_n}]$ is pinched.
Similarly, $\gamma_{-\infty} :=\lim_{n\ra\infty}\gamma_{-n}$ defines
an u.s.c.\ invariant graph which is pinched to $\varphi_{A_{n-1}}^-$.

Let $\theta^*$ be a fiber on which $\tilde{A}$ is pinched.
Due to the compactness of $\tilde{A}$ we can choose %
$\delta_2 \in (0,\delta_1)$ for which
\begin{equation}
\label{eq:pinchedball2}
  [\gamma_\infty(\theta),\varphi^+_{A_n}(\theta)] \ \ssq \
  B_{\eps/8}(\varphi^+_{A_n}(\theta^*)) \quad \forall \theta
  \in B_{\delta_2}(\theta^*) \ .
\end{equation}
Fix $k_1 \in \N$ such that for some $\theta' \in
B_{\delta_2/2}(\theta_0)$ there holds $\gamma_{-k_1}(\theta') \leq
x_0$.  Note that such $k_1$ and $\theta'$ exist even if $x_0 <
\gamma_{-\infty}(\theta_0)$, because arbitrarily close to $\theta_0$ 
there are fibers where $[\varphi_{A_{n-1}}^-,\gamma_{-\infty}]$ is 
pinched.  Further, fix $K$ and $\delta_3 \in (0,\delta_2/2)$ with the
property that
\begin{equation} \label{eq:prop2}
   \gamma_j(\theta) \ \geq \
\varphi^+_{A_n}(\theta^*) - \eps/4 \quad \forall j \geq K,\ \theta\in
B_{\delta_3(\theta^*)} \ .
\end{equation}
 Choose $k_2 \geq K$ such that
$\theta'+(k_1+k_2)\omega \in B_{\delta_3}(\theta^*)$ and let $k:=
k_1+k_2$. Then
\begin{eqnarray*}
    F_{\eps,\theta_0}^k(x_0) & \stackrel{(\ref{eq:bound2})}{\geq} &
    F_{0,\theta'}^k(x_0) +  \frac{\eps}{2} \ \stackrel{}{\geq} \
    F^k_{0,\theta'}(\gamma_{-k_1}(\theta')) + \frac{\eps}{2} \ 
    = \ \gamma_{k_2}(\theta'+k\omega) + \frac{\eps}{2} \\
    & \stackrel{(\ref{eq:prop2})}{\geq} & \varphi^+_{A_n}(\theta^*) +
    \frac{\eps}{4} \  \stackrel{~(\ref{eq:pinchedball2})}{>} \
    \varphi^+_{A_n}(\theta_0+k\omega)\ . 
  \end{eqnarray*}

\ \\
\underline{Case $({\cal A}2)$:} \ \emph{$A_n$ is defined by
(\ref{eq:phin+1}) and (\ref{eq:An+1}).} \bigskip

By definition (\ref{eq:An+1}) we have $\varphi_{A_n}^+ = \psi^+$, and
hence $\widehat{A}:=[\psi,\varphi_{A_n}^+]$ is pinched by
Lemma~\ref{lem:pinched} (recall that $\psi$ is l.s.c.\ by definition).
As before, we choose some fiber $\theta^*$ on which $\widehat{A}$ is
pinched and fix $\delta_2\in (0,\delta_1)$ such that
\begin{equation}
\label{eq:pinchedball3}
  [\psi(\theta),\varphi^+_{A_n}(\theta)] \ \ssq \
  B_{\eps/8}(\varphi^+_{A_n}(\theta^*)) \quad \forall \theta
  \in B_{\delta_2}(\theta^*) \ .
\end{equation}
Further, we choose $\theta' \in \kreis$ which satisfies
$\varphi^+_{A_{n-1}}(\theta') \geq \psi(\theta')-\eps/2$
(see~(\ref{eq:theta'})). By Lemma~\ref{lem:reflexivegraphs} we have
$\varphi^+_{A_{n-1}} = \varphi^{-+}_{A_{n-1}}$ and therefore
$(\theta',\varphi^+_{A_{n-1}}(\theta')) \in
\overline{\Phi^-_{A_{n-1}}}$.  Consequently, due to the l.s.c.\ of
$\varphi^-_{A_{n-1}}$, there exists an open set $U \ssq
B_{\delta_2/2}(\theta')$ such that
\begin{equation}
  \label{eq:prop3}
  \varphi^+_{A_{n-1}}(\theta) \ \geq \ \varphi^-_{A_{n-1}}(\theta) \
  \geq \ \psi(\theta') - \eps \quad \forall \theta \in U \ .
\end{equation}
Choose $k_1 \geq 1$ for which $\theta_0 + k_1\omega \in U$. Then
\begin{equation} \label{eq:prop4}
F_{\eps,\theta_0}^{k_1}(x_0) \ \geq \ F_{0,\theta_0}^{k_1}(x_0)+\eps \ 
\geq \ \varphi^+_{A_{n-1}}(\theta_0+k_1\omega) + \eps \ 
\stackrel{~(\ref{eq:prop3})}{\geq} \ \psi(\theta') \ .
\end{equation}
Now let $k_2 \geq 1$ such that $\theta'+k_2 \omega \in
B_{\delta_2/2}(\theta^*)$ and define $k:= k_1+k_2$. Then
\begin{eqnarray*}
F^k_{\eps,\theta_0}(x_0) & = &
F^{k_2}_{\eps,\theta_0+k_1\omega}(F_{\eps,\theta_0}^{k_1}(x_0)) \ 
\stackrel{~(\ref{eq:prop4})}{\geq} \
F_{\eps,\theta_0+k_1\omega}^{k_2}(\psi(\theta')) \\ 
& \stackrel{~(\ref{eq:bound2})}{\geq} &
F_{0,\theta'}^{k_2}(\psi(\theta')) + \frac{\eps}{2} \ = \
\psi(\theta'+k_2\omega) + \frac{\eps}{2}\
\stackrel{~(\ref{eq:pinchedball3})}{>} \
\varphi^+_{A_n}(\theta_0+k\omega) \ .
\end{eqnarray*}

\qed


\end{document}